\documentclass[12pt]{article}

\usepackage{amsfonts,amssymb,amsthm,amsmath,latexsym}
\usepackage{subeqnarray,indent,url,cite,bm}
\usepackage{tensor}  

\date{March 29, 2019 \\[1.5mm] revised September 24, 2019}

\begin{document}

\title{\vspace*{-1cm}
       Wall's continued-fraction characterization \\
       of Hausdorff moment sequences: \\[3mm]
       A conceptual proof
      }

\author{
     {\small Alan D.~Sokal}                  \\[2mm]
     {\small\it Department of Mathematics}   \\[-2mm]
     {\small\it University College London}   \\[-2mm]
     {\small\it Gower Street}                \\[-2mm]
     {\small\it London WC1E 6BT}             \\[-2mm]
     {\small\it UNITED KINGDOM}              \\[-2mm]
     {\small\tt sokal@math.ucl.ac.uk}        \\[-2mm]
     {\protect\makebox[5in]{\quad}}  
     \\[-2mm]
     {\small\it Department of Physics}       \\[-2mm]
     {\small\it New York University}         \\[-2mm]
     {\small\it 4 Washington Place}          \\[-2mm]
     {\small\it New York, NY 10003}          \\[-2mm]
     {\small\it USA}                         \\[-2mm]
     {\small\tt sokal@nyu.edu}               \\[3mm]
}

\maketitle
\thispagestyle{empty}   

\vspace*{-5mm}

\begin{abstract}
I give an elementary proof of Wall's continued-fraction characterization
of Hausdorff moment sequences.
\end{abstract}


\bigskip
\bigskip
\noindent
{\bf Mathematics Subject Classification (MSC 2010) codes:}
44A60 (Primary);  05A15, 30B70, 30E05 (Secondary).

\clearpage

\newtheorem{theorem}{Theorem}[section]
\newtheorem{proposition}[theorem]{Proposition}
\newtheorem{lemma}[theorem]{Lemma}
\newtheorem{corollary}[theorem]{Corollary}
\newtheorem{definition}[theorem]{Definition}
\newtheorem{conjecture}[theorem]{Conjecture}
\newtheorem{question}[theorem]{Question}
\newtheorem{problem}[theorem]{Problem}
\newtheorem{example}[theorem]{Example}

\renewcommand{\theenumi}{\alph{enumi}}
\renewcommand{\labelenumi}{(\theenumi)}
\def\eop{\hbox{\kern1pt\vrule height6pt width4pt
depth1pt\kern1pt}\medskip}
\def\prf{\par\noindent{\bf Proof.\enspace}\rm}
\def\rmk{\par\medskip\noindent{\bf Remark\enspace}\rm}

\newcommand{\textbfit}[1]{\textbf{\textit{#1}}}

\newcommand{\bigdash}{%
\smallskip\begin{center} \rule{5cm}{0.1mm} \end{center}\smallskip}

\newcommand{\safepar}{ {\protect\hfill\protect\break\hspace*{5mm}} }

\newcommand{\be}{\begin{equation}}
\newcommand{\ee}{\end{equation}}
\newcommand{\<}{\langle}
\renewcommand{\>}{\rangle}
\newcommand{\widebar}{\overline}
\def\reff#1{(\protect\ref{#1})}
\def\spose#1{\hbox to 0pt{#1\hss}}
\def\ltapprox{\mathrel{\spose{\lower 3pt\hbox{$\mathchar"218$}}
    \raise 2.0pt\hbox{$\mathchar"13C$}}}
\def\gtapprox{\mathrel{\spose{\lower 3pt\hbox{$\mathchar"218$}}
    \raise 2.0pt\hbox{$\mathchar"13E$}}}
\def\textprime{${}^\prime$}
\def\proof{\par\medskip\noindent{\sc Proof.\ }}
\def\firstproof{\par\medskip\noindent{\sc First Proof.\ }}
\def\secondproof{\par\medskip\noindent{\sc Second Proof.\ }}
\def\alternateproof{\par\medskip\noindent{\sc Alternate Proof.\ }}
\def\algebraicproof{\par\medskip\noindent{\sc Algebraic Proof.\ }}
\def\combinatorialproof{\par\medskip\noindent{\sc Combinatorial Proof.\ }}
\def\proofof#1{\bigskip\noindent{\sc Proof of #1.\ }}
\def\firstproofof#1{\bigskip\noindent{\sc First Proof of #1.\ }}
\def\secondproofof#1{\bigskip\noindent{\sc Second Proof of #1.\ }}
\def\thirdproofof#1{\bigskip\noindent{\sc Third Proof of #1.\ }}
\def\algebraicproofof#1{\bigskip\noindent{\sc Algebraic Proof of #1.\ }}
\def\combinatorialproofof#1{\bigskip\noindent{\sc Combinatorial Proof of #1.\ }}
\def\sketchofproof{\par\medskip\noindent{\sc Sketch of proof.\ }}
\renewcommand{\qed}{ $\square$ \bigskip}
\newcommand{\myendremark}{ $\blacksquare$ \bigskip}
\def\half{ {1 \over 2} }
\def\third{ {1 \over 3} }
\def\twothird{ {2 \over 3} }
\def\smfrac#1#2{{\textstyle{#1\over #2}}}
\def\smhalf{ {\smfrac{1}{2}} }
\def\smquarter{ {\smfrac{1}{4}} }
\newcommand{\real}{\mathop{\rm Re}\nolimits}
\renewcommand{\Re}{\mathop{\rm Re}\nolimits}
\newcommand{\imag}{\mathop{\rm Im}\nolimits}
\renewcommand{\Im}{\mathop{\rm Im}\nolimits}
\newcommand{\sgn}{\mathop{\rm sgn}\nolimits}
\newcommand{\tr}{\mathop{\rm tr}\nolimits}
\newcommand{\tg}{\mathop{\rm tg}\nolimits}
\newcommand{\supp}{\mathop{\rm supp}\nolimits}
\newcommand{\disc}{\mathop{\rm disc}\nolimits}
\newcommand{\diag}{\mathop{\rm diag}\nolimits}
\newcommand{\csch}{\mathop{\rm csch}\nolimits}
\newcommand{\tridiag}{\mathop{\rm tridiag}\nolimits}
\newcommand{\AZ}{\mathop{\rm AZ}\nolimits}
\newcommand{\perm}{\mathop{\rm perm}\nolimits}
\def\hboxscript#1{ {\hbox{\scriptsize\em #1}} }
\renewcommand{\emptyset}{\varnothing}
\newcommand{\eqdef}{\stackrel{\rm def}{=}}

\newcommand{\restrict}{\upharpoonright}

\newcommand{\compinv}{{\langle -1 \rangle}}   

\newcommand{\scra}{{\mathcal{A}}}
\newcommand{\scrb}{{\mathcal{B}}}
\newcommand{\scrc}{{\mathcal{C}}}
\newcommand{\scrd}{{\mathcal{D}}}
\newcommand{\scre}{{\mathcal{E}}}
\newcommand{\scrf}{{\mathcal{F}}}
\newcommand{\scrg}{{\mathcal{G}}}
\newcommand{\scrh}{{\mathcal{H}}}
\newcommand{\scri}{{\mathcal{I}}}
\newcommand{\scrj}{{\mathcal{J}}}
\newcommand{\scrk}{{\mathcal{K}}}
\newcommand{\scrl}{{\mathcal{L}}}
\newcommand{\scrm}{{\mathcal{M}}}
\newcommand{\scrn}{{\mathcal{N}}}
\newcommand{\scro}{{\mathcal{O}}}
\newcommand{\scrp}{{\mathcal{P}}}
\newcommand{\scrq}{{\mathcal{Q}}}
\newcommand{\scrr}{{\mathcal{R}}}
\newcommand{\scrs}{{\mathcal{S}}}
\newcommand{\scrt}{{\mathcal{T}}}
\newcommand{\scrv}{{\mathcal{V}}}
\newcommand{\scrw}{{\mathcal{W}}}
\newcommand{\scrz}{{\mathcal{Z}}}

\newcommand{\bfa}{{\mathbf{a}}}
\newcommand{\bfb}{{\mathbf{b}}}
\newcommand{\bfc}{{\mathbf{c}}}
\newcommand{\bfd}{{\mathbf{d}}}
\newcommand{\bfe}{{\mathbf{e}}}
\newcommand{\bfj}{{\mathbf{j}}}
\newcommand{\bfi}{{\mathbf{i}}}
\newcommand{\bfk}{{\mathbf{k}}}
\newcommand{\bfl}{{\mathbf{l}}}
\newcommand{\bfm}{{\mathbf{m}}}
\newcommand{\bfx}{{\mathbf{x}}}
\renewcommand{\k}{{\mathbf{k}}}
\newcommand{\n}{{\mathbf{n}}}
\newcommand{\vv}{{\mathbf{v}}}
\newcommand{\bv}{{\mathbf{v}}}
\newcommand{\w}{{\mathbf{w}}}
\newcommand{\x}{{\mathbf{x}}}
\newcommand{\cc}{{\mathbf{c}}}
\newcommand{\zero}{{\mathbf{0}}}
\newcommand{\one}{{\mathbf{1}}}
\newcommand{\bmm}{{\mathbf{m}}}

\newcommand{\ahat}{{\widehat{a}}}
\newcommand{\Zhat}{{\widehat{Z}}}

\newcommand{\C}{{\mathbb C}}
\newcommand{\D}{{\mathbb D}}
\newcommand{\Z}{{\mathbb Z}}
\newcommand{\N}{{\mathbb N}}
\newcommand{\Q}{{\mathbb Q}}
\newcommand{\PP}{{\mathbb P}}
\newcommand{\R}{{\mathbb R}}
\newcommand{\RR}{{\mathbb R}}
\newcommand{\E}{{\mathbb E}}

\newcommand{\Sym}{{\mathfrak{S}}}
\newcommand{\SymB}{{\mathfrak{B}}}
\newcommand{\Alt}{{\mathrm{Alt}}}

\newcommand{\myle}{\preceq}
\newcommand{\myge}{\succeq}
\newcommand{\mygt}{\succ}

\newcommand{\B}{{\sf B}}
\newcommand{\OB}{{\sf OB}}
\newcommand{\OS}{{\sf OS}}
\newcommand{\OO}{{\sf O}}
\newcommand{\SP}{{\sf SP}}
\newcommand{\OSP}{{\sf OSP}}
\newcommand{\Eu}{{\sf Eu}}
\newcommand{\ERR}{{\sf ERR}}
\newcommand{\sfB}{{\sf B}}
\newcommand{\sfD}{{\sf D}}
\newcommand{\sfE}{{\sf E}}
\newcommand{\sfG}{{\sf G}}
\newcommand{\sfJ}{{\sf J}}
\newcommand{\sfP}{{\sf P}}
\newcommand{\sfQ}{{\sf Q}}
\newcommand{\sfS}{{\sf S}}
\newcommand{\sfT}{{\sf T}}
\newcommand{\sfW}{{\sf W}}
\newcommand{\sfMV}{{\sf MV}}
\newcommand{\AMV}{{\sf AMV}}
\newcommand{\BM}{{\sf BM}}
\newcommand{\NC}{{\sf NC}}

\newcommand{\emIB}{{\hbox{\em IB}}}
\newcommand{\emIP}{{\hbox{\em IP}}}
\newcommand{\emOB}{{\hbox{\em OB}}}
\newcommand{\emSC}{{\hbox{\em SC}}}

\newcommand{\stat}{{\rm stat}}
\newcommand{\cs}{{\rm cs}}
\newcommand{\cyc}{{\rm cyc}}
\newcommand{\Asc}{{\rm Asc}}
\newcommand{\asc}{{\rm asc}}
\newcommand{\Des}{{\rm Des}}
\newcommand{\des}{{\rm des}}
\newcommand{\Exc}{{\rm Exc}}
\newcommand{\exc}{{\rm exc}}
\newcommand{\Wex}{{\rm Wex}}
\newcommand{\wex}{{\rm wex}}
\newcommand{\Fix}{{\rm Fix}}
\newcommand{\fix}{{\rm fix}}
\newcommand{\lrmax}{{\rm lrmax}}
\newcommand{\rlmax}{{\rm rlmax}}
\newcommand{\Rec}{{\rm Rec}}
\newcommand{\rec}{{\rm rec}}
\newcommand{\Arec}{{\rm Arec}}
\newcommand{\arec}{{\rm arec}}
\newcommand{\ERec}{{\rm ERec}}
\newcommand{\erec}{{\rm erec}}
\newcommand{\EArec}{{\rm EArec}}
\newcommand{\earec}{{\rm earec}}
\newcommand{\recarec}{{\rm recarec}}
\newcommand{\nonrec}{{\rm nonrec}}
\newcommand{\Cpeak}{{\rm Cpeak}}
\newcommand{\cpeak}{{\rm cpeak}}
\newcommand{\Cval}{{\rm Cval}}
\newcommand{\cval}{{\rm cval}}
\newcommand{\Cdasc}{{\rm Cdasc}}
\newcommand{\cdasc}{{\rm cdasc}}
\newcommand{\Cddes}{{\rm Cddes}}
\newcommand{\cddes}{{\rm cddes}}
\newcommand{\cdrise}{{\rm cdrise}}
\newcommand{\cdfall}{{\rm cdfall}}
\newcommand{\Peak}{{\rm Peak}}
\newcommand{\peak}{{\rm peak}}
\newcommand{\Val}{{\rm Val}}
\newcommand{\val}{{\rm val}}
\newcommand{\Dasc}{{\rm Dasc}}
\newcommand{\dasc}{{\rm dasc}}
\newcommand{\Ddes}{{\rm Ddes}}
\newcommand{\ddes}{{\rm ddes}}
\newcommand{\inv}{{\rm inv}}
\newcommand{\maj}{{\rm maj}}
\newcommand{\rs}{{\rm rs}}
\newcommand{\cross}{{\rm cr}}
\newcommand{\crosshat}{{\widehat{\rm cr}}}
\newcommand{\nest}{{\rm ne}}
\newcommand{\rodd}{{\rm rodd}}
\newcommand{\reven}{{\rm reven}}
\newcommand{\lodd}{{\rm lodd}}
\newcommand{\leven}{{\rm leven}}
\newcommand{\sg}{{\rm sg}}
\newcommand{\bl}{{\rm bl}}
\newcommand{\tran}{{\rm tr}}
\newcommand{\area}{{\rm area}}
\newcommand{\ret}{{\rm ret}}
\newcommand{\peaks}{{\rm peaks}}
\newcommand{\hl}{{\rm hl}}
\newcommand{\sll}{{\rm sl}}
\newcommand{\negg}{{\rm neg}}
\newcommand{\imp}{{\rm imp}}
\newcommand{\osg}{{\rm osg}}
\newcommand{\ons}{{\rm ons}}
\newcommand{\isg}{{\rm isg}}
\newcommand{\ins}{{\rm ins}}
\newcommand{\LL}{{\rm LL}}
\newcommand{\height}{{\rm ht}}
\newcommand{\as}{{\rm as}}

\newcommand{\ba}{{\bm{a}}}
\newcommand{\bahat}{{\widehat{\bm{a}}}}
\newcommand{\sfa}{{{\sf a}}}
\newcommand{\bb}{{\bm{b}}}
\newcommand{\bc}{{\bm{c}}}
\newcommand{\bchat}{{\widehat{\bm{c}}}}
\newcommand{\bd}{{\bm{d}}}
\newcommand{\bee}{{\bm{e}}}
\newcommand{\bff}{{\bm{f}}}
\newcommand{\bg}{{\bm{g}}}
\newcommand{\bh}{{\bm{h}}}
\newcommand{\bll}{{\bm{\ell}}}
\newcommand{\bp}{{\bm{p}}}
\newcommand{\br}{{\bm{r}}}
\newcommand{\bs}{{\bm{s}}}
\newcommand{\bu}{{\bm{u}}}
\newcommand{\bw}{{\bm{w}}}
\newcommand{\bx}{{\bm{x}}}
\newcommand{\by}{{\bm{y}}}
\newcommand{\bz}{{\bm{z}}}
\newcommand{\bA}{{\bm{A}}}
\newcommand{\bB}{{\bm{B}}}
\newcommand{\bC}{{\bm{C}}}
\newcommand{\bE}{{\bm{E}}}
\newcommand{\bF}{{\bm{F}}}
\newcommand{\bG}{{\bm{G}}}
\newcommand{\bH}{{\bm{H}}}
\newcommand{\bI}{{\bm{I}}}
\newcommand{\bJ}{{\bm{J}}}
\newcommand{\bM}{{\bm{M}}}
\newcommand{\bN}{{\bm{N}}}
\newcommand{\bP}{{\bm{P}}}
\newcommand{\bQ}{{\bm{Q}}}
\newcommand{\bS}{{\bm{S}}}
\newcommand{\bT}{{\bm{T}}}
\newcommand{\bW}{{\bm{W}}}
\newcommand{\bX}{{\bm{X}}}
\newcommand{\bIB}{{\bm{IB}}}
\newcommand{\bOB}{{\bm{OB}}}
\newcommand{\bOS}{{\bm{OS}}}
\newcommand{\bERR}{{\bm{ERR}}}
\newcommand{\bSP}{{\bm{SP}}}
\newcommand{\bMV}{{\bm{MV}}}
\newcommand{\bBM}{{\bm{BM}}}
\newcommand{\balpha}{{\bm{\alpha}}}
\newcommand{\bbeta}{{\bm{\beta}}}
\newcommand{\bgamma}{{\bm{\gamma}}}
\newcommand{\bdelta}{{\bm{\delta}}}
\newcommand{\bkappa}{{\bm{\kappa}}}
\newcommand{\bomega}{{\bm{\omega}}}
\newcommand{\bsigma}{{\bm{\sigma}}}
\newcommand{\btau}{{\bm{\tau}}}
\newcommand{\bpsi}{{\bm{\psi}}}
\newcommand{\bzeta}{{\bm{\zeta}}}
\newcommand{\bone}{{\bm{1}}}
\newcommand{\bzero}{{\bm{0}}}

\newcommand{\Cbar}{{\overline{C}}}
\newcommand{\Dbar}{{\overline{D}}}
\newcommand{\dbar}{{\overline{d}}}
\def\Ctilde{{\widetilde{C}}}
\def\Etilde{{\widetilde{E}}}
\def\Ftilde{{\widetilde{F}}}
\def\Gtilde{{\widetilde{G}}}
\def\Htilde{{\widetilde{H}}}
\def\Ptilde{{\widetilde{P}}}
\def\Chat{{\widehat{C}}}
\def\ctilde{{\widetilde{c}}}
\def\zbar{{\overline{Z}}}
\def\pitilde{{\widetilde{\pi}}}

\newcommand{\sech}{{\rm sech}}

%
%
\newcommand{\sn}{{\rm sn}}
\newcommand{\cn}{{\rm cn}}
\newcommand{\dn}{{\rm dn}}
\newcommand{\sm}{{\rm sm}}
\newcommand{\cm}{{\rm cm}}

%
%
\newcommand{\zfz}{ {{}_0 \! F_0} }
\newcommand{\zfo}{ {{}_0  F_1} }
\newcommand{\ofz}{ {{}_1 \! F_0} }
\newcommand{\ofo}{ {{}_1 \! F_1} }
\newcommand{\oft}{ {{}_1 \! F_2} }

%
%
\newcommand{\FHyper}[2]{ {\tensor[_{#1 \!}]{F}{_{#2}}\!} }
\newcommand{\FHYPER}[5]{ {\FHyper{#1}{#2} \!\biggl(
   \!\!\begin{array}{c} #3 \\[1mm] #4 \end{array}\! \bigg|\, #5 \! \biggr)} }
\newcommand{\tfo}{ {\FHyper{2}{1}} }
\newcommand{\tfz}{ {\FHyper{2}{0}} }
\newcommand{\threefz}{ {\FHyper{3}{0}} }
\newcommand{\FHYPERbottomzero}[3]{ {\FHyper{#1}{0} \!\biggl(
   \!\!\begin{array}{c} #2 \\[1mm] \hbox{---} \end{array}\! \bigg|\, #3 \! \biggr)} }
\newcommand{\FHYPERtopzero}[3]{ {\FHyper{0}{#1} \!\biggl(
   \!\!\begin{array}{c} \hbox{---} \\[1mm] #2 \end{array}\! \bigg|\, #3 \! \biggr)} }

\newcommand{\phiHyper}[2]{ {\tensor[_{#1}]{\phi}{_{#2}}} }
\newcommand{\psiHyper}[2]{ {\tensor[_{#1}]{\psi}{_{#2}}} }
\newcommand{\PhiHyper}[2]{ {\tensor[_{#1}]{\Phi}{_{#2}}} }
\newcommand{\PsiHyper}[2]{ {\tensor[_{#1}]{\Psi}{_{#2}}} }
\newcommand{\phiHYPER}[6]{ {\phiHyper{#1}{#2} \!\left(
   \!\!\begin{array}{c} #3 \\ #4 \end{array}\! ;\, #5, \, #6 \! \right)\!} }
\newcommand{\psiHYPER}[6]{ {\psiHyper{#1}{#2} \!\left(
   \!\!\begin{array}{c} #3 \\ #4 \end{array}\! ;\, #5, \, #6 \! \right)} }
\newcommand{\PhiHYPER}[5]{ {\PhiHyper{#1}{#2} \!\left(
   \!\!\begin{array}{c} #3 \\ #4 \end{array}\! ;\, #5 \! \right)\!} }
\newcommand{\PsiHYPER}[5]{ {\PsiHyper{#1}{#2} \!\left(
   \!\!\begin{array}{c} #3 \\ #4 \end{array}\! ;\, #5 \! \right)\!} }
\newcommand{\zerophizero}{ {\phiHyper{0}{0}} }
\newcommand{\ophizero}{ {\phiHyper{1}{0}} }
\newcommand{\zphio}{ {\phiHyper{0}{1}} }
\newcommand{\ophio}{ {\phiHyper{1}{1}} }
\newcommand{\tphio}{ {\phiHyper{2}{1}} }
\newcommand{\tphiz}{ {\phiHyper{2}{0}} }
\newcommand{\tPhio}{ {\PhiHyper{2}{1}} }
\newcommand{\opsio}{ {\psiHyper{1}{1}} }

%
%
\newcommand{\stirlingsubset}[2]{\genfrac{\{}{\}}{0pt}{}{#1}{#2}}
\newcommand{\stirlingcycleold}[2]{\genfrac{[}{]}{0pt}{}{#1}{#2}}
\newcommand{\stirlingcycle}[2]{\left[\! \stirlingcycleold{#1}{#2} \!\right]}
\newcommand{\assocstirlingsubset}[3]{{\genfrac{\{}{\}}{0pt}{}{#1}{#2}}_{\! \ge #3}}
\newcommand{\genstirlingsubset}[4]{{\genfrac{\{}{\}}{0pt}{}{#1}{#2}}_{\! #3,#4}}
\newcommand{\euler}[2]{\genfrac{\langle}{\rangle}{0pt}{}{#1}{#2}}
\newcommand{\eulergen}[3]{{\genfrac{\langle}{\rangle}{0pt}{}{#1}{#2}}_{\! #3}}
\newcommand{\eulersecond}[2]{\left\langle\!\! \euler{#1}{#2} \!\!\right\rangle}
\newcommand{\eulersecondgen}[3]{{\left\langle\!\! \euler{#1}{#2} \!\!\right\rangle}_{\! #3}}
\newcommand{\binomvert}[2]{\genfrac{\vert}{\vert}{0pt}{}{#1}{#2}}
\newcommand{\binomsquare}[2]{\genfrac{[}{]}{0pt}{}{#1}{#2}}


\newenvironment{sarray}{
             \textfont0=\scriptfont0
             \scriptfont0=\scriptscriptfont0
             \textfont1=\scriptfont1
             \scriptfont1=\scriptscriptfont1
             \textfont2=\scriptfont2
             \scriptfont2=\scriptscriptfont2
             \textfont3=\scriptfont3
             \scriptfont3=\scriptscriptfont3
           \renewcommand{\arraystretch}{0.7}
           \begin{array}{l}}{\end{array}}

\newenvironment{scarray}{
             \textfont0=\scriptfont0
             \scriptfont0=\scriptscriptfont0
             \textfont1=\scriptfont1
             \scriptfont1=\scriptscriptfont1
             \textfont2=\scriptfont2
             \scriptfont2=\scriptscriptfont2
             \textfont3=\scriptfont3
             \scriptfont3=\scriptscriptfont3
           \renewcommand{\arraystretch}{0.7}
           \begin{array}{c}}{\end{array}}


\newcommand*\circled[1]{\tikz[baseline=(char.base)]{
  \node[shape=circle,draw,inner sep=1pt] (char) {#1};}}
\newcommand{\ostar}{{\circledast}}
\newcommand{\ostarN}{{\,\circledast_{\vphantom{\dot{N}}N}\,}}
\newcommand{\ostarPsi}{{\,\circledast_{\vphantom{\dot{\Psi}}\Psi}\,}}
\newcommand{\starN}{{\,\ast_{\vphantom{\dot{N}}N}\,}}
\newcommand{\starpsi}{{\,\ast_{\vphantom{\dot{\bpsi}}\!\bpsi}\,}}
\newcommand{\starone}{{\,\ast_{\vphantom{\dot{1}}1}\,}}
\newcommand{\startwo}{{\,\ast_{\vphantom{\dot{2}}2}\,}}
\newcommand{\starinfty}{{\,\ast_{\vphantom{\dot{\infty}}\infty}\,}}
\newcommand{\starT}{{\,\ast_{\vphantom{\dot{T}}T}\,}}

\clearpage

Let us recall that a sequence $\ba = (a_n)_{n \ge 0}$ of real numbers
is called a Hamburger (resp.\ Stieltjes, resp.\ Hausdorff) moment sequence
\cite{Shohat_43,Akhiezer_65,Berg_84,Simon_98,Schmudgen_17}
if there exists a positive measure $\mu$
on $\R$ (resp.\ on $[0,\infty)$, resp.\ on $[0,1]$)
such that $a_n = \int \! x^n \, d\mu(x)$ for all $n \ge 0$.
%
One fundamental characterization of Stieltjes moment sequences
was found by Stieltjes \cite{Stieltjes_1894} in 1894
(see also \cite[pp.~327--329]{Wall_48}):
A sequence $\ba = (a_n)_{n \ge 0}$ of real numbers
is a Stieltjes moment sequence if and only if
there exist real numbers $\alpha_0,\alpha_1,\alpha_2,\ldots \ge 0$
such that
\be
   \sum_{n=0}^{\infty} a_n t^n
   \;=\;
   \cfrac{\alpha_0}{1 - \cfrac{\alpha_1 t}{1 - \cfrac{\alpha_2 t}{1 - \cdots}}}
   \label{eq.Sfrac}
\ee
in the sense of formal power series.
(That is, the ordinary generating function
 $f(t) = \sum\limits_{n=0}^{\infty} a_n t^n$
 can be represented as a Stieltjes-type continued fraction
 with nonnegative coefficients.)
Moreover, the coefficients $\balpha = (\alpha_i)_{i \ge 0}$ are unique
if we make the convention that $\alpha_i = 0$ implies
$\alpha_j = 0$ for all $j > i$;
we shall call such a sequence $\balpha$ {\em standard}\/.

Since every Hausdorff moment sequence is a Stieltjes moment sequence,
its ordinary generating function clearly has a continued-fraction expansion
of the form \reff{eq.Sfrac} with coefficients $\balpha \ge 0$.
But {\em which}\/ sequences $\balpha \ge 0$
correspond to Hausdorff moment sequences?
The answer was given by Wall \cite[Theorems~4.1 and 6.1]{Wall_40} in 1940:
A~sequence $\ba = (a_n)_{n \ge 0}$ of real numbers
is a Hausdorff moment sequence if and only if
there exist real numbers $c \ge 0$ and $g_1,g_2,g_3,\ldots \in [0,1]$
such that
\be
   \sum_{n=0}^{\infty} a_n t^n
   \;=\;
   \cfrac{c}{1 - \cfrac{g_1 t}{1 - \cfrac{(1-g_1)g_2 t}{1 -  \cfrac{(1-g_2)g_3 t}{1- \cfrac{(1-g_3)g_4 t}{1-\cdots}}}}}
   \label{eq.thm.wall}
\ee
in the sense of formal power series.

Wall's proof of this result was based on an interesting
but somewhat mysterious identity for continued fractions
\cite[Theorem~2.1]{Wall_40}
together with some complex-analysis arguments.\footnote{
   See also \cite{Kim_00} for a combinatorial proof of Wall's identity,
   and see \cite{Randrianarivony_94b,Zeng_96,Han_99}
   for some interesting applications of it.
}
Four years later, Wall \cite{Wall_44} gave a new proof,
based on Schur's \cite{Schur_17} characterization of
analytic functions bounded in the unit disc
and the Herglotz--Riesz \cite{Herglotz_11,Riesz_11} integral representation
of analytic functions in the unit disc with positive real part.

Here I would like to present an alternate proof of Wall's theorem
that is not only very simple
but also gives insight into why the coefficients in \reff{eq.thm.wall}
take the form $\alpha_n = (1 - g_{n-1}) g_n$.

This proof requires two well-known elementary facts about moment sequences:

1) $\ba$ is a Stieltjes moment sequence if and only if
the ``aerated'' sequence $\widehat{\ba} = (a_0,0,a_1,0,a_2,0,\ldots)$
is a Hamburger moment sequence.
Indeed, the even subsequence of a Hamburger moment sequence
is always a Stieltjes moment sequence;
and conversely, if $\ba$ is a Stieltjes moment sequence
that is represented by a measure $\mu$ supported on $[0,\infty)$,
then $\widehat{\ba}$ is represented by the even measure
$\widehat{\mu} = (\tau^+ + \tau^-)/2$ on $\R$,
where $\tau^\pm$ is the image of $\mu$
under the map $x \mapsto \pm \sqrt{x}$.
In particular, if $\mu$ is supported on $[0,\Lambda]$,
then $\widehat{\mu}$ is supported on $[-\sqrt{\Lambda},\sqrt{\Lambda}]$.

2) If the Hamburger moment sequence $\ba = (a_n)_{n \ge 0}$
satisfies $|a_n| \le A B^n$ with $A,B < \infty$,
then the representing measure $\mu$ is unique
and is supported on $[-B,B]$.
In particular, a Hausdorff moment sequence
always has a unique representing measure.
(In fact, the representing measure $\mu$ of a Hamburger moment sequence is unique
under the vastly weaker hypothesis $|a_n| \le A B^n n!$,
or even under the yet weaker hypothesis
$\sum\limits_{n=1}^\infty a_{2n}^{-1/2n} = \infty$
\cite[section~4.2]{Schmudgen_17};
but we shall not need these latter results.)

Besides Stieltjes-type continued fractions \reff{eq.Sfrac}
[henceforth called S-fractions for short],
we shall also make use of Jacobi-type continued fractions (J-fractions)
\be
   f(t)
   \;=\;
   \cfrac{1}{1 - \gamma_0 t - \cfrac{\beta_1 t^2}{1 - \gamma_1 t - \cfrac{\beta_2 t^2}{1 - \cdots}}}
  \label{eq.Jfrac}
\ee
(always considered as formal power series in the indeterminate $t$).\footnote{
   My use of the terms ``S-fraction'' and ``J-fraction'' follows
   the general practice in the combinatorial literature,
   starting with Flajolet \cite{Flajolet_80}.
   The classical literature on continued fractions
   \cite{Perron,Wall_48,Jones_80,Lorentzen_92,Cuyt_08}
   generally uses a different terminology.
   For instance, Jones and Thron \cite[pp.~128--129, 386--389]{Jones_80}
   use the term ``regular C-fraction''
   for (a minor variant of) what I~have called an S-fraction,
   and the term ``associated continued fraction''
   for (a minor variant of) what I~have called a J-fraction.
}
We shall need three elementary facts about these continued fractions:

1) The contraction formula:  We have
\be
   \cfrac{1}{1 - \cfrac{\alpha_1 t}{1 - \cfrac{\alpha_2 t}{1 -  \cfrac{\alpha_3 t}{1- \cdots}}}}
   \;\;=\;\;
   \cfrac{1}{1 - \alpha_1 t - \cfrac{\alpha_1 \alpha_2 t^2}{1 - (\alpha_2 + \alpha_3) t - \cfrac{\alpha_3 \alpha_4 t^2}{1 - (\alpha_4 + \alpha_5) t - \cfrac{\alpha_5 \alpha_6 t^2}{1- \cdots}}}}
 \label{eq.contraction_even}
\ee
as an identity between formal power series.
In other words, an S-fraction with coefficients $\balpha$
is equal to a J-fraction with coefficients $\bgamma$ and $\bbeta$, where
\begin{subeqnarray}
   \gamma_0  & = &  \alpha_1
       \slabel{eq.contraction_even.coeffs.a}   \\
   \gamma_n  & = &  \alpha_{2n} + \alpha_{2n+1}  \quad\hbox{for $n \ge 1$}
       \slabel{eq.contraction_even.coeffs.b}   \\
   \beta_n  & = &  \alpha_{2n-1} \alpha_{2n}
       \slabel{eq.contraction_even.coeffs.c}
 \label{eq.contraction_even.coeffs}
\end{subeqnarray}
See \cite[pp.~20--22]{Wall_48} for the classic algebraic proof of
the contraction formula \reff{eq.contraction_even};
see \cite[Lemmas~1 and 2]{Dumont_94b}
\cite[proof of Lemma~1]{Dumont_95} \cite[Lemma~4.5]{DiFrancesco_10}
for a very simple variant algebraic proof;
and see \cite[pp.~V-31--V-32]{Viennot_83}
for an enlightening combinatorial proof,
based on Flajolet's \cite{Flajolet_80} combinatorial interpretation
of S-fractions (resp.\ J-fractions)
as generating functions for Dyck (resp.\ Motzkin) paths
with height-dependent weights.

2) Binomial transform:  Fix a real number $\xi$,
and let $\ba = (a_n)_{n \ge 0}$ be a sequence of real numbers.
Then the $\xi$-binomial transform of $\ba$
is defined to be the sequence $\bb = (b_n)_{n \ge 0}$ given by
\be
   b_n  \;=\; \sum_{k=0}^n \binom{n}{k} \, a_k \, \xi^{n-k}
   \;.
\ee
Note that if $a_n = \int \! x^n \, d\mu(x)$,
then $b_n = \int (x + \xi)^n \, d\mu(x)$.
In other words, if $\ba$ is a Hamburger moment sequence
with representing measure $\mu$,
then $\bb$ is a Hamburger moment sequence
with representing measure $T_\xi \mu$
(the $\xi$-translate of $\mu$).

Now suppose that the ordinary generating function of $\ba$
is given by a J-fraction:
\be
   \sum_{n=0}^{\infty} a_n t^n
   \;=\;
   \cfrac{1}{1 - \gamma_0 t - \cfrac{\beta_1 t^2}{1 - \gamma_1 t - \cfrac{\beta_2 t^2}{1 - \cdots}}}
   \;.
  \label{eq.Jfrac.bis}
\ee
Then the $\xi$-binomial transform $\bb$ of $\ba$ is given by a J-fraction
in which we make the replacement $\gamma_i \to \gamma_i + \xi$:
\be
   \sum_{n=0}^{\infty} b_n t^n
   \;=\;
   \cfrac{1}{1 - (\gamma_0 + \xi) t - \cfrac{\beta_1 t^2}{1 - (\gamma_1 + \xi) t - \cfrac{\beta_2 t^2}{1 - \cdots}}}
   \;.
  \label{eq.Jfrac.bis2}
\ee
See \cite[Proposition~4]{Barry_09} for an algebraic proof
of \reff{eq.Jfrac.bis2};
or see \cite{Sokal_totalpos} for a simple combinatorial proof
based on Flajolet's \cite{Flajolet_80} theory.

3) An upper bound:  If $\ba$ is given by the S-fraction \reff{eq.Sfrac}
with $0 \le \alpha_i \le 1$ for all~$i$,
then $0 \le a_n \le C_n \le 4^n$,
where $C_n = \displaystyle \frac{1}{n+1} \binom{2n}{n}$
is the $n$th Catalan number.
The proof is simple:
If we consider the coefficients $\balpha$ in \reff{eq.Sfrac}
to be indeterminates,
then it is easy to see that $a_n$ is a polynomial in $\alpha_0,\ldots,\alpha_n$
with nonnegative integer coefficients
(namely, $\alpha_0$ times the Stieltjes--Rogers polynomial
 $S_n(\balpha)$ \cite{Flajolet_80});
so $a_n$ is an increasing function of $\balpha$ on the set $\balpha \ge 0$.
On the other hand, if $\alpha_i = 1$ for all $i$,
then \reff{eq.Sfrac} represents a series $f(t) = \sum_{n=0}^\infty a_n t^n$
satisfying $f(t) = 1/[1 - t f(t)]$,
from which it follows that $f(t) = [1 - \sqrt{1-4t}]/(2t)$
and hence (by binomial expansion) that $a_n = C_n$.

\bigskip

\proofof{Wall's theorem}
Let $\ba = (a_n)_{n \ge 0}$ be a Hausdorff moment sequence;
we can assume without loss of generality that $a_0 = 1$.
Then $\ba$ has a (unique) representing measure $\mu$ supported on $[0,1]$,
and its ordinary generating function is given by a unique S-fraction
\reff{eq.Sfrac} with $\alpha_0 = 1$ and standard coefficients $\balpha \ge 0$.
Now let $\widehat{\ba} = (a_0,0,a_1,0,a_2,0,\ldots)$
be the aerated sequence;  it is a Hamburger moment sequence
with a (unique) even representing measure $\widehat{\mu}$ supported on $[-1,1]$,
and its ordinary generating function is given by the J-fraction
with coefficients $\bgamma = 0$ and $\bbeta = \balpha$:
\be
   \sum_{n=0}^{\infty} \widehat{a}_n t^n
   \;=\;
   \sum_{n=0}^{\infty} a_n t^{2n}
   \;=\;
   \cfrac{1}{1 - \cfrac{\alpha_1 t^2}{1 - \cfrac{\alpha_2 t^2}{1 - \cdots}}}
   \;.
\ee
Now let $\widetilde{\ba}$ be the 1-binomial transform of $\widehat{\ba}$;
it is a Stieltjes moment sequence
with a (unique) representing measure $\widetilde{\mu} = T_1 \widehat{\mu}$
supported on $[0,2]$,
and its ordinary generating function is given by a J-fraction
with coefficients $\bgamma = \bone$ and $\bbeta = \balpha$:
\be
   \sum_{n=0}^{\infty} \widetilde{a}_n t^n
   \;=\;
   \cfrac{1}{1 - t - \cfrac{\alpha_1 t^2}{1 - t - \cfrac{\alpha_2 t^2}{1 - \cdots}}}
   \;.
\ee
But since $\widetilde{\ba}$ is a Stieltjes moment sequence,
its ordinary generating function is also given by an S-fraction
with nonnegative coefficients, call them $\balpha'$.
Comparing the J-fraction and the S-fraction using the contraction formula
\reff{eq.contraction_even}/\reff{eq.contraction_even.coeffs}, we see that
\begin{subeqnarray}
   & &
   1 \,=\, \alpha'_1 \,=\, \alpha'_2 + \alpha'_3 \,=\, \alpha'_4 + \alpha'_5
       \,=\, \ldots
      \slabel{eq.alphaprime.contraction.a}   \\[2mm]
   & &
   \alpha_1 \,=\, \alpha'_1 \alpha'_2 \,,\quad
   \alpha_2 \,=\, \alpha'_3 \alpha'_4 \,,\quad \ldots
\end{subeqnarray}
It follows from \reff{eq.alphaprime.contraction.a}
that $\alpha'_i \in [0,1]$ for all $i$.
Setting $g_n = \alpha'_{2n}$
shows that $\alpha_1 = g_1$ and $\alpha_n = (1-g_{n-1}) g_n$ for $n \ge 2$,
which is precisely the representation \reff{eq.thm.wall}.

Conversely, suppose that $\ba$ is given by an S-fraction \reff{eq.thm.wall}
with coefficients $c=1$ and $g_i \in [0,1]$.
Then $\ba$ is a Stieltjes moment sequence satisfying $a_n \le 4^n$,
so that the representing measure $\mu$ is unique and is supported on $[0,4]$.
(Of course, we will soon see that $\mu$ is actually supported on $[0,1]$.)
Then the aerated sequence $\widehat{\ba} = (a_0,0,a_1,0,a_2,0,\ldots)$
is a Hamburger moment sequence
with a unique representing measure $\widehat{\mu}$
that is even and supported on $[-2,2]$,
and its ordinary generating function is given by a J-fraction
with coefficients $\bgamma = 0$, $\beta_1 = g_1$
and $\beta_n = (1-g_{n-1}) g_n$ for $n \ge 2$.
Now let $\widetilde{\ba}$ be the 1-binomial transform of $\widehat{\ba}$:
it is a Hamburger moment sequence with a
unique representing measure $\widetilde{\mu} = T_1 \widehat{\mu}$
supported on $[-1,3]$,
and its ordinary generating function is given by a J-fraction with coefficients
$\bgamma = \bone$, $\beta_1 = g_1$
and $\beta_n = (1-g_{n-1}) g_n$ for $n \ge 2$.
But the contraction formula
\reff{eq.contraction_even}/\reff{eq.contraction_even.coeffs}
shows that this J-fraction is equivalent to an S-fraction with coefficients
$\alpha'_1 = 1$ and $\alpha'_{2n} = g_n$, $\alpha'_{2n+1} = 1 - g_n$
for $n \ge 1$.
Since all these coefficients are nonnegative,
it follows that $\widetilde{\ba}$ is a Stieltjes moment sequence.
Therefore $\widetilde{\mu}$ is supported on $[0,3]$,
so that $\widehat{\mu} = T_{-1} \widetilde{\mu}$ is supported on $[-1,2]$.
But since $\widehat{\mu}$ is even, it must actually be supported on $[-1,1]$.
Hence $\mu$ is supported on $[0,1]$,
which shows that $\ba$ is a Hausdorff moment sequence.
\qed

\section*{Acknowledgments}

I wish to thank Alexander Dyachenko, Lily Liu, Mathias P\'etr\'eolle
and Jiang Zeng for helpful conversations.
This research was supported in part by
Engineering and Physical Sciences Research Council grant EP/N025636/1.

\end{document}